\documentclass[A4paper,12pt]{article}

 \usepackage{latexsym,epsfig,amsthm, amscd, a4wide}

 \usepackage{amsfonts,amsmath,mathrsfs,amssymb}
 \usepackage{times,helvet,courier,type1cm}
 \usepackage{color, graphicx}
 \usepackage[colorlinks, linkcolor=blue, citecolor=blue]{hyperref}
 \usepackage{hyperref}

 \allowdisplaybreaks

\allowdisplaybreaks[4]

 \newtheorem{thm0}{Theorem}[section]

 \newtheorem{def1}[thm0]{Definition}
 \newtheorem{lem1}[thm0]{Lemma}
 \newtheorem{thm1}[thm0]{Theorem}
 \newtheorem{cor1}[thm0]{Corollary}
 \newtheorem{pro1}[thm0]{Proposition}
 \newtheorem{con1}[thm0]{Condition}
 \newtheorem{rem1}[thm0]{Remark}

 \def\blemma{\begin{lem1}}\def\elemma{\end{lem1}}
 \def\btheorem{\begin{thm1}}\def\etheorem{\end{thm1}}
 \def\bcorollary{\begin{cor1}}\def\ecorollary{\end{cor1}}
 \def\bproposition{\begin{pro1}}\def\eproposition{\end{pro1}}
 
 \def\bremark{\begin{rem1}}\def\eremark{\end{rem1}}

 \def\bproof{\begin{proof}~}\def\eproof{\hfill$\square$\end{proof}}

 \def\itDelta{{\it\Delta}}

 \def\itOmega{{\mathit\Omega}}

 \def\benumerate{\begin{enumerate}}\def\eenumerate{\end{enumerate}}
 \def\bitemize{\begin{itemize}}\def\eitemize{\end{itemize}}
 \def\bibitm{\bibitem}

 \def\beqlb{\begin{eqnarray}}\def\eeqlb{\end{eqnarray}}
 \def\beqnn{\begin{eqnarray*}}\def\eeqnn{\end{eqnarray*}}

 \DeclareMathOperator{\sgn}{sgn}

 \def\eqref#1{{\rm(\ref{#1})}}

 \def\aar{&}
 \def\qqquad{\qquad\qquad}

 \def\mbf{\mathbf}\def\mrm{\mathrm}\def\mfr{\mathfrak}
 \def\mbb{\mathbb}\def\mcr{\mathscr}

 \def\qed{\hfill$\square$\smallskip}

 \def\<{\langle}\def\>{\rangle}
 \def\vepsilon{\varepsilon}

 \def\d{\mrm{d}}\def\e{\mrm{e}}\def\im{\mrm{i}}

 \def\supp{\mrm{supp}}

 \def\qqquad{\qquad\quad}

\begin{document}

\noindent{(Version: 2025-11-25)}

\bigskip\bigskip

\noindent\textbf{\LARGE Stochastic integral representations for the}

\smallskip

\noindent\textbf{\LARGE Ray--Knight theorem of the L\'evy forest}

\bigskip

\noindent{Pei-Sen Li$^{\,a)}$, Zenghu Li$^{\,b)}$ and Wenjing Zhang$^{\,b)}$}

\smallskip\noindent{a)~\it School of Mathematics and Statistics, Beijing Institute of Technology, \\ Beijing 100872, China}

\smallskip\noindent{b)~\it Laboratory of Mathematics and Complex Systems, \\ School of Mathematical Sciences, Beijing Normal University, \\ Beijing 100875, China}

\smallskip\noindent{Emails: peisenli@bit.edu.cn, lizh@bnu.edu.cn, zhangwenjing@mail.bnu.edu.cn}

\bigskip\bigskip

{\narrower

\noindent\textbf{Abstract:~} We present a simple stochastic integral representation for the local times of the height process of a spectrally positive L\'evy process stopped at a hitting time. From the representation we derive a strong stochastic equation for the local time process of the type of Bertoin and Le~Gall (Illinois J. Math., 2006) and Dawson and Li (Ann. Probab., 2012). This leads to a representation of the Ray--Knight theorem of Le~Gall and Le~Jan (Ann. Probab., 1998) and Duquesne and Le~Gall (Ast\'erisque, 2002), which codes the genealogical forest of a continuous-state branching process. The results extend those in the recent work of A\"id\'ekon et al.\ (Sci. China Math., 2024) for a Brownian motion with a local time drift.

\bigskip

\noindent\textbf{Keywords and phrases:~} Spectrally positive L\'evy process, height process, local time, stochastic integral representation, Ray--Knight theorem, continuous-state branching process, stochastic equation, genealogical forest, Tanaka formula.

\medskip

\noindent\textbf{2020 Mathematics Subject Classification:~} 60G51, 60J55, 60J80

}

\bigskip

\section{Introduction}

\setcounter{equation}{0}

A \textit{continuous-state branching process} (CB-process) models the stochastic evolution of a large population with small individuals and arises naturally as the scaling limit of discrete Galton--Watson branching processes. Suppose that $\alpha\ge 0$ and $\beta\ge 0$ are real constants and $(z\land z^2)\pi(\d z)$ is a finite measure on $(0,\infty)$. Let $\psi$ be the positive ($=\,$nonnegative) function on $\mbb{R}_+:= [0,\infty)$ defined by
 \begin{align}\label{psi(lam)=def}
\psi(\lambda)= \alpha\lambda + \beta\lambda^2+\int_0^\infty (\e^{-\lambda z}-1+\lambda z) \pi(\d z), \quad \lambda\ge 0.
 \end{align}
By a CB-process with (subcritical) \textit{branching mechanism} $\psi$ we mean a positive Markov process with Feller transition semigroup $(Q_t)_{t\ge 0}$ defined by
 \begin{align}\label{LaplaceQ_t(x,.)}
\int_{[0,\infty)} \e^{-\lambda y} Q_t(x,\d y)
 =
\e^{-xv_t(\lambda)}, \qquad \lambda\ge 0,
 \end{align}
where $t\mapsto v_t(\lambda)$ is the unique positive solution of
 \begin{align*}%\label{cumulanteq}
\frac{\d v_t}{\d t}(\lambda)= -\psi(v_t(\lambda)), \quad v_0(\lambda)= \lambda.
 \end{align*}
In the special case where $\pi(0,\infty)= 0$, the process is known as a \textit{Feller branching diffusion}, which was first studied by Feller \cite{Fel51}.

Suppose that $W(\d s,\d u)$ is a time-space Gaussian white noise on $(0,\infty)^2$ with intensity $\d s\d u$ and $\tilde{M}(\d s,\d z,\d u)$ is a compensated time-space Poisson random measure on $(0,\infty)^3$ with intensity $\d s\pi(\d z)\d u$. Let $L(\d s,\d u)$ be the time-space L\'evy noise on $(0,\infty)^2$ defined by
 \begin{align}\label{L(ds,du)=def}
L(\d s,\d u)= -\alpha\d s\d u + \sqrt{2\beta}W(\d s,\d u) + \int_{\{0< z< \infty\}} z \tilde{M}(\d s,\d z,\d u).
 \end{align}
By Dawson and Li \cite[Theorem~3.1]{DaL12}, for any $x\ge 0$ we can construct a CB-process $\{X_t(x): t\ge 0\}$ with branching mechanism $\psi$ by the pathwise unique solution to the stochastic equation
 \begin{align}\label{DLEX_t(x)=x+..}
X_t(x)= x + \int_0^t\int_0^{X_{s-}(x)} L(\d s,\d u), \quad t\ge 0;
 \end{align}
see also Bertoin and Le~Gall \cite[Proposition~2]{BeL06}. Here and in the sequel, we understand that
 \begin{align*}
\int_a^b= \int_{(a,b]}, ~~ \int_a^\infty= \int_{(a,\infty)}, \quad a\le b\in \mbb{R}.
 \end{align*}
According to \cite[Theorem~3.3]{DaL12}, the path-valued process $x\mapsto \{X_t(x): t\ge 0\}$ has stationary and independent increments. The solution flow $\{X_t(x): x,t\ge 0\}$ of \eqref{DLEX_t(x)=x+..} captures the genealogical structures of the population represented by the CB-process. Stochastic flows with similar properties were used by Bertoin and Le~Gall \cite{BeL00, BeL03, BeL05, BeL06} to study the coalescent process introduced by Bolthausen and Sznitman \cite{BoS98} in their work on Ruelle's continuous probability cascades.

A method of coding the genealogical forest of the CB-process was suggested by Le~Gall and Le~Jan \cite{LeL98} and Duquesne and Le~Gall \cite{DLG02}. They represented the forest by the \textit{height process} of a spectrally positive L\'evy process with Laplace exponent given by \eqref{psi(lam)=def} and proved Ray--Knight theorems for the local times of the height process. The Ray--Knight theorem states that a weak solution flow to \eqref{DLEX_t(x)=x+..} is given by the local times of the height process. Their approach provides powerful tools in the research and has led to many deep results. Similar structures of related models have been studied by many authors; see, e.g., \cite{BFF18, Lam02, LPW13, LPW22, Par16, Xu24, Xu24+}. In particular, Berestycki et al.\ \cite{BFF18} proposed a novel genealogical description for CB-processes with competition based on interactive pruning of L\'evy trees, and established a Ray--Knight theorem for the processes in terms of the local times of suitably pruned forests. A different description for the genealogical forests of CB-processes with interaction was given by Li et al.\ \cite{LPW22}. In the special case where $\beta>0$, they proved a stochastic integral representation for the height process and gave a characterization for its local times by a Tanaka formula. The corresponding Ray--Knight theorem was proved in \cite{LPW22} by considering an approximating sequence of processes with finite L\'evy measures.

A very nice representation for the classical Ray--Knight theorem was given recently by A\"id\'ekon et al.\ \cite{AHS24}. Starting from a Brownian motion with a local time drift, they constructed explicitly a space-time Gaussian white noise by a stochastic integral relative to the Brownian motion so that the local times of the stopped motion satisfy stochastic equations of the type of \eqref{DLEX_t(x)=x+..}. In fact, they showed that the stochastic equations are reformulations of the classical Tanaka formulae. The representation given in \cite{AHS24} has played the key role in the exploration of the Brownian loop soup on the real line in \cite{AHS25+} by the same authors.

In this work, we are interested in stochastic integral representations of the Ray--Knight theorem in the setting of Le~Gall and Le~Jan \cite{LeL98} and Duquesne and Le~Gall \cite{DLG02}. Let $(\itOmega, \mcr{G}, \mbf{P})$ be a complete probability space furnished with the filtration $(\mcr{G}_t)_{t\ge 0}$ satisfying the usual hypotheses. Suppose that $(B_t: t\ge 0)$ is a standard $(\mcr{G}_t)$-Brownian motion and $N(\d s,\d z)$ is a time-space $(\mcr{G}_t)$-Poisson random measure on $(0,\infty)^2$ with intensity $\d s\pi(\d z)$. Let $(\xi_t: t\ge 0)$ be the spectrally positive $(\mcr{G}_t)$-L\'evy process defined by the \textit{L\'evy--It\^o representation}:
 \begin{align}\label{xi=(Levy-Ito)}
\xi_t= - \alpha t + \sqrt{2\beta}B_t + \int_0^t\int_0^\infty z\tilde{N}(\d s,\d z), \quad t\ge 0,
 \end{align}
where $\tilde{N}(\d s,\d z)$ denotes the compensated random measure. It is well-known that
 \begin{align*}%\label{Laplacexi}
\mbf{E}[\e^{-\lambda(\xi_{s+t}-\xi_s)}]
 =
\e^{t\psi(\lambda)}, \quad s\ge 0, t\ge 0, \lambda\ge 0.
 \end{align*}
The L\'evy process given by \eqref{xi=(Levy-Ito)} has a c\`adl\`ag version. We shall always use this version of the process. By the construction of the Poisson random measure, there is a random countable subset $\{(t_i,z_i): i\in I\}$ of $(0,\infty)^2$, where each $t_i$ is a $(\mcr{G}_t)$-stopping time, such that
 \begin{align}\label{N=sum_I}
N(\d s,\d z)= \sum_{i\in I} \delta_{(t_i,z_i)}(\d s,\d z), \quad s,z>0.
 \end{align}
Then $(\xi_t: t\ge 0)$ has discontinuity set $J:= \{t_i: i\in I\}$ and it makes a jump of size $z_i> 0$ at time $t_i\in J$. We refer to the books \cite{Ber96, Kyp14} for the theory of L\'evy processes.

\smallskip\noindent\textbf{Definition.~} Suppose that $(\xi_t: t\ge 0)$ is a real c\`adl\`ag stochastic process. Let $C(\mbb{R})$ be set of bounded continuous functions on $\mbb{R}$. A two-parameter process $\{L_t(a,\xi): a\in \mbb{R}, t\ge 0\}$ is called the \textit{local time} of $(\xi_t: t\ge 0)$ if for every $t\ge 0$ and $f\in C(\mbb{R})$ we have a.s.
 \begin{align}\label{loctimeeq}
\int_\mbb{R} L_t(x,\xi)f(x)\d x
 =
\int_0^t f(\xi_s)\d s.
 \end{align}

Let $(\xi_t: t\ge 0)$ be the spectrally positive L\'evy process defined by \eqref{xi=(Levy-Ito)}. Consider the \textit{supremum process} $(S_t: t\ge 0)$ and the \textit{reflected process} $(R_t: t\ge 0)$ defined by
 \begin{align}\label{defS_tR_t}
S_t= \sup_{0\le s\le t} \xi_s,
 \quad
R_t= S_t - \xi_t.
 \end{align}
It is known that the positive semi-martingale $(R_t: t\ge 0)$ has a jointly measurable local time $(t,a)\mapsto L_t(a,R)$ that is continuous in $t\ge 0$; see, e.g., \cite[pp.212-213 and p.216]{Pro04}. For any $t\ge 0$ define the \textit{time-reversal}
 \begin{align}\label{hat{xi}_t^{(t)}=def}
\hat{\xi}_s^{(t)}= \xi_t-\xi_{(t-s)-},
 \quad
0\le s\le t,
 \end{align}
where we understand $\xi_{0-}= 0$ by convention. Then $(\hat{\xi}_s^{(t)}: 0\le s\le t)$ is distributed identically with $(\xi_s: 0\le s\le t)$. Let $(\hat{S}_s^{(t)}: 0\le s\le t)$ and $(\hat{R}_s^{(t)}: 0\le s\le t)$ be the supremum process and the reflected process of $(\hat{\xi}_s^{(t)}: 0\le s\le t)$, respectively. Following Duquesne and Le~Gall \cite[Definition~1.2.1 and Lemma~1.2.1]{DLG02}, we define the \textit{height process} $(H_t: t\ge 0)$ by
 \begin{align}\label{H_t=L_t(0,hR^(t))}
H_t= L_t(0,\hat{R}^{(t)}), \quad t\ge 0.
 \end{align}

By \cite[Theorem~1.4.3]{DLG02}, the height process has a continuous modification if and only if the following \textit{Grey's condition} is satisfied:
 \begin{align}\label{Grey'scondition}
\int_1^{\infty} \frac{\d u}{\psi(u)}< \infty.
 \end{align}
\textit{In this work, we always assume the above Grey's condition is satisfied and consider a continuous version of the height process.}

As pointed out by Duquesne and Le~Gall \cite[p.7, p.31]{DLG02}, the height process $(H_t: t\ge 0)$ is neither a semimartingale nor a Markov process unless $\pi(0,\infty)= 0$. However, by \cite[Proposition~1.3.3]{DLG02}, it still has a jointly measurable local time $\{L_t(a): a,t\ge 0\}$ such that $t\mapsto L_t(a)$ is continuous and increasing. Let $\mcr{L}_2$ be the set of Borel functions $f$ on $(0,\infty)^2$ such that
 \begin{align}\label{iintf<infty}
\int_0^t \d s\int_0^\infty f(s,u)^2\d u< \infty, \quad t\ge 0.
 \end{align}
Let $\mcr{L}_3(\pi)$ be the set of Borel functions $g$ on $(0,\infty)^3$ such that
 \begin{align}\label{iiintg<infty}
\int_0^t \d s\int_0^\infty \pi(\d z)\int_0^\infty g(s,z,u)^2\d u< \infty, \quad t\ge 0.
 \end{align}
For $x\ge 0$ define the \textit{hitting time} $T_x= \inf\{t\ge 0: \xi_t= -x\}$. The main results of the paper are the following:

\btheorem\label{thL_{T_x}(a)=x+int..} The local time process $\{L_{T_x}(t): t\ge 0\}$ has the stochastic integral representation
 \begin{align}\label{L_{T_x}(t)=x+int_0^{T_x}}
L_{T_x}(t)= x + \int_0^{T_x} 1_{\{H_s\le t\}}\d \xi_s, \quad t\ge 0.
 \end{align}
\etheorem

\bproposition\label{th-GP-noises} There is a time-space Gaussian white noise $W(\d s,\d u)$ on $(0,\infty)^2$ with intensity $\d s\d u$ and a compensated time-space Poisson random measure $\tilde{M}(\d s,\d z,\d u)$ on $(0,\infty)^3$ with intensity $\d s\pi(\d z)\d u$ such that, for $t\ge 0$, $f\in \mcr{L}_2$ and $g\in \mcr{L}_3(\pi)$,
 \begin{align}\label{iintfdW=def}
\int_0^t\int_0^\infty f(s,u) W(\d s,\d u)
 =
\int_0^\infty 1_{\{H_s\le t\}} f(H_s,L_s(H_s))\d B_s
 \end{align}
and
 \begin{align}\label{iintgdM=def}
\aar\int_0^t\int_0^\infty\int_0^\infty g(s,z,u) \tilde{M}(\d s,\d z,\d u) \cr
 \aar\qqquad\qqquad
= \int_0^\infty\int_0^\infty 1_{\{H_s\le t\}} g(H_s,z,L_s(H_s)) \tilde{N}(\d s,\d z).
 \end{align}
\eproposition

\btheorem\label{th-DL-inteq} Let $L(\d s,\d u)$ be the time-space L\'evy noise defined by \eqref{L(ds,du)=def}, where $W(\d s,\d u)$ and $\tilde{M}(\d s,\d z,\d u)$ are given by \eqref{iintfdW=def} and \eqref{iintgdM=def}, respectively. Then the process $t\mapsto X_t(x):= L_{T_x}(t)$ has a c\`adl\`ag modification, which is the pathwise unique solution to \eqref{DLEX_t(x)=x+..}. \etheorem

Theorems~\ref{thL_{T_x}(a)=x+int..} and~\ref{th-DL-inteq} give stochastic integral representations for the stochastic flow $\{L_{T_x}(t): x\ge 0, t\ge 0\}$. For $\beta> 0$, we prove \eqref{L_{T_x}(t)=x+int_0^{T_x}} by a Tanaka formula for the local time of the height process, which was essentially established by Li et al.\ \cite{LPW22}. For $\beta= 0$, the proof of \eqref{L_{T_x}(t)=x+int_0^{T_x}} is based on a semimartingale characterization the CB-process $\{L_{T_x}(t): t\ge 0\}$ given in \cite{Li20}. The following example illustrates the connections of the results presented above with those obtained by A\"id\'ekon et al.\ \cite{AHS24}.

\noindent\textbf{Example.~} Let $(B_t: t\ge 0)$ be a Brownian motion and let $\mfr{L}_t= L_t(0,B)$ be its local time at level zero accumulated by time $t\ge 0$. Clearly, we also have $\mfr{L}_t= L_t(0,2|B|)$ for every $t\ge 0$. By the classical Tanaka formula, we can define another Brownian motion $(\xi_t: t\ge 0)$ by
 \begin{align*}
\xi_t= \int_0^t \sgn(B_s)\d B_s= |B_t| - \mfr{L}_t.
 \end{align*}
Then $(\xi_t: t\ge 0)$ has height process $(2|B_t|: t\ge 0)$ by Theorem~\ref{th-Hproc-inteq}. According to Proposition~\ref{th-GP-noises}, there is a time-space Gaussian white noise $W(\d s,\d u)$ on $(0,\infty)^2$ with intensity $\d s\d u$ such that, for any $f\in \mcr{L}_2$,
 \begin{align}\label{iintfdW=intfdB}
\int_0^t\int_0^\infty f(s,u)W(\d s,\d u)
 =
\int_0^\infty 1_{\{2|B_s|\le t\}}f(2|B_s|,L_s(2|B_s|,2|B|)) \sgn(B_s)\d B_s.
 \end{align}
By Theorem~\ref{th-DL-inteq}, the local time process  $\{L_{T_x}(t,2|B|): t\ge 0\}$ solves the equation
 \begin{align*}
L_{T_x}(t,2|B|)= x + \int_0^t\int_0^{L_{T_x}(s,2|B|)} W(\d s,\d u), \quad t\ge 0.
 \end{align*}
This is essentially a spacial form of the stochastic equation given by A\"id\'ekon et al.\ \cite[Theorem~1.2-(i)]{AHS24}, who studied the more general \textit{$\mu$-process} $(\xi_t^\mu: t\ge 0)$ for a real constant $\mu\neq 0$ defined by
 \begin{align*}
\xi_t^\mu= |B_t| - \mu\mfr{L}_t, \quad t\ge 0.
 \end{align*}
The relation \eqref{iintfdW=intfdB} is an insightful discovery of A\"id\'ekon et al.\ \cite{AHS24}.

The paper is organized as follows. In Section~2, we prove some stochastic equations for the height process and its local time. The proof of Theorem~\ref{thL_{T_x}(a)=x+int..} is given in Section~3. In Section~4, we introduce some $\sigma$-algebras indexed by the spatial parameter. The measurabilities of some important random variables relative to those $\sigma$-algebras are discussed. The proofs of Proposition~\ref{th-GP-noises} and Theorem~\ref{th-DL-inteq} are given in Sections~5 and~6, respectively.

\section{Preliminaries}

\setcounter{equation}{0}

Throughout this section, we assume $\beta> 0$. In this case, Li et al.\ \cite{LPW22} proved a stochastic integral representation for the height process and a Tanaka formula for its local times. We shall give new treatments of those results and derive some useful stochastic equations from them.

\bproposition\label{th-betL_t(0,R)=..} The supremum process $t\mapsto S_t$ has continuous part $t\mapsto \beta L_t(0,R)$, that is,
 \begin{align}\label{betL_t(0,R)=..}
\beta L_t(0,R)= S_t - \sum_{0< t_i\le t} \Delta S_{t_i}, \quad t\ge 0,
 \end{align}
where
 \begin{align}\label{DelS_{t_i}=def}
\Delta S_{t_i}:= S_{t_i}-S_{t_i-}= (\xi_{t_i}-S_{t_i-})^+.
 \end{align}
\eproposition

\bproof Observe that $t\mapsto S_t$ is an increasing c\`adl\`ag process making a possible jump at time $t_i\in J$ with jump size $\Delta S_{t_i}$ given by \eqref{DelS_{t_i}=def}. Let $S_t^c$ and $R_t^c$ denote the continuous part of $S_t$ and $R_t$, respectively. Then
 \begin{align*}%\label{R_t=S_t_c+..}
R_t^c= S_t^c + \alpha t - \sqrt{2\beta} B_t,
 \quad
R_t= R_t^c + \sum_{t_i\le t} \Delta R_{t_i},
 \end{align*}
where $\Delta R_{t_i}= \Delta S_{t_i} - \Delta \xi_{t_i}\le 0$. Let $\{\mcr{L}_t(a): t,a\ge 0\}$ be the semimartingale local time of $(R_t: t\ge 0)$. For any $t\ge 0$ and $f\in C[0,\infty)$, we have
 \begin{align*}
\int_0^\infty \mcr{L}_t(x)f(x)\d x
 =
2\beta\int_0^t f(R_s)\d s;
 \end{align*}
see, e.g., \cite[pp.213--214 and p.216]{Pro04}. Then $L_t(a,R)= (2\beta)^{-1}\mcr{L}_t(a)$ defines a local time of $(R_t: t\ge 0)$ in the scale specified in \eqref{loctimeeq}. Since $t\mapsto R_t$ can only have negative jumps, an application of Tanaka's formula shows that
 \begin{align*}%\label{R_t=int1_.dR+..}
\beta L_t(0,R)= R_t - \int_0^t 1_{\{R_{s-}> 0\}}\d R_s
 =
\int_0^t 1_{\{R_{s-}= 0\}}\d R_s
 =
\int_0^t 1_{\{R_{s-}= 0\}}\d R_s^c;
 \end{align*}
see, e.g., \cite[p.213]{Pro04}. Since $1_{\{R_t= 0\}}= 0$ for a.e.\ $t\ge 0$, we obtain
 \begin{align}\label{betL(0,R)=int1dR}
\beta L_t(0,R)= \int_0^t 1_{\{R_{s-}= 0\}}\d (S_s^c + \alpha s - \sqrt{2\beta} B_s)
 =
\int_0^t 1_{\{R_{s-}= 0\}}\d S_s^c.
 \end{align}
It is easy to show that $G:= \{t >0: R_{t-}> 0$ and $R_t> 0\}$ is an open subset of $(0,\infty)$, so it is the union of a sequence of disjoint open intervals, say $(l_k,r_k)$, $k= 1,2, \cdots$. For any $t\in (l_k,r_k)$ we have $R_t= S_t - \xi_t> 0$, and so $t\mapsto S_t$ remains constant in $(l_k,r_k)$. We also use $S^c$ to denote the locally finite measure on $[0,\infty)$ such that $S^c[0,t]= S^c(t)$ for $t\ge 0$. Then $S^c(l_k,r_k)= 0$, and hence $S^c(G)= 0$. Recall that $J= \{t_i: i\in I\}$ denotes the countable discontinuity set of $(\xi_t: t\ge 0)$. The continuity of $t\mapsto S_t$ implies that $S^c(J)= 0$. Since $\{t >0: R_{t-}> 0\}\subset G\cup J$, we have $S^c(\{t >0: R_{t-}> 0\})= 0$ by the subadditivity. Then $\beta L_t(0,R)= S_t^c$ by \eqref{betL(0,R)=int1dR}. \eproof

We believe that the result of Proposition~\ref{th-betL_t(0,R)=..} is known to the expert in L\'evy processes, but we could not find it in the form of \eqref{betL_t(0,R)=..} in the literature. The above proof, included for the convenience of the reader, identifies the right coefficient $\beta$ on the l.h.s.\ of \eqref{betL_t(0,R)=..}. See \cite{Ber96, Kyp14} for much more intensive discussions on the subject for spectrally negative L\'evy processes.

The next theorem was proved in Li et al.\ \cite{LPW22} by an approximation of $(\xi_t: t\ge 0)$ using a sequence of processes with finite L\'evy measures. Here we provide a simpler proof of the result based on Proposition~\ref{th-betL_t(0,R)=..}.

\btheorem\label{th-Hproc} The height process $(H_t: t\ge 0)$ defined by \eqref{H_t=L_t(0,hR^(t))} has the stochastic integral representation
 \begin{align}\label{Hproc-intN}
\beta H_t= \xi_t - I_0(t) - \int_0^t\int_0^\infty \big(z+I_s(t)-\xi_s\big)^+ N(\d s,\d z), \quad t\ge 0,
 \end{align}
where $I_s(t)= \inf_{s\le u\le t} \xi_u$. \etheorem

\bproof By applying Proposition~\ref{th-betL_t(0,R)=..} to the spectrally positive L\'evy process $(\hat{\xi}_s^{(t)}: 0\le s\le t)$ we have
 \begin{align}\label{betL_t(0,hat{R}^{(t)})}
\beta H_t= \beta L_t(0,\hat{R}^{(t)})= \hat{S}_t^{(t)} - \sum_{0< t_i\le t} \Delta \hat{S}_{t-t_i}^{(t)},
 \end{align}
where
 \begin{align*}
\hat{S}_t^{(t)}= \sup_{0\le s\le t}(\xi_t - \xi_{(t-s)-})
 \end{align*}
and
 \begin{align*}
\Delta \hat{S}_{t-t_i}^{(t)}
 =
\big(\hat{\xi}_{t-t_i}^{(t)} - \hat{S}_{(t-t_i)-}^{(t)}\big)^+.
 \end{align*}
By elementary calculations,
 \begin{align*}
\hat{S}_t^{(t)}
 =
\xi_t - \inf_{0\le s\le t} \xi_{(t-s)-}
 =
\xi_t - I_0(t).
 \end{align*}
In view of \eqref{DelS_{t_i}=def}, we have
 \begin{align*}
\Delta \hat{S}_{t-t_i}^{(t)}
 \aar=
\Big(\inf_{0\le s< t-t_i}\xi_{(t-s)-} - \xi_{t_i-}\Big)^+ \cr
 \aar=
\big(I_{t_i}(t) - \xi_{t_i-}\big)^+
 =
\big(z_i + I_{t_i}(t) - \xi_{t_i}\big)^+.
 \end{align*}
Then \eqref{Hproc-intN} follows by \eqref{betL_t(0,hat{R}^{(t)})}. \eproof

Clearly, the path $s\mapsto I_s(t)$ increases from $I_0(t)$ to $I_t(t)= \xi_t$ on the interval $[0,t]$. It is continuous on the set $[0,t]\setminus J$ and makes a possible jump of size $\itDelta I_{t_i}(t):= I_{t_i}(t) - \xi_{t_i-}(t)$ at time $t_i\in [0,t]\cap J$. Observe that
 \begin{align}\label{DelI_{t_i}(t)}
\itDelta I_{t_i}(t)= \big(I_{t_i}(t) - \xi_{t_i-}\big)^+= \big(z_i + I_{t_i}(t) - \xi_{t_i}\big)^+.
 \end{align}
Then we can rewrite \eqref{Hproc-intN} as
 \begin{align}\label{Hproc-sum0}
\beta H_t= \xi_t - I_0(t) - \sum_{0< t_i\le t} \itDelta I_{t_i}(t).
 \end{align}
It follows that
 \begin{align}\label{betH=leb(..)}
\beta H_t= \mfr{m}(\{I_s(t): 0\le s\le t\}),
 \end{align}
where $\mfr{m}$ denotes the Lebesgue measure; see also \cite[Formula (14)]{DLG02}.

\bremark\label{thHproc-sum1} A continuous version of the height process is given by \eqref{betH=leb(..)}; see \cite[p.24]{DLG02}. For $i\in I$ define the stopping time $r(t_i)= \inf\{t\ge t_i: \xi_t - \xi_{t_i}= -z_i\}$. Then we can rewrite \eqref{Hproc-intN} or \eqref{Hproc-sum0} as
 \begin{align}\label{Hproc-sum2}
\beta H_t= \xi_t - I_0(t) - \sum_{0< t_i\le t} \big(z_i + I_{t_i}(t\land r(t_i)) - \xi_{t_i}\big),
 \end{align}
where $t\mapsto - I_0(t)$ reflects the L\'evy process above level zero and $t\mapsto \xi_{t_i} - I_{t_i}(t\land r(t_i))$ reflects the process in the interval $(t_i,r(t_i))$ above level $\xi_{t_i}$. For any given $t_i< s< r(t_i)$, the conditional law of $(H_{t\land r(t_i)}: t\ge s)$ depends not only on $H_s$ but also on $H_{t_i}$. Therefore, the height process is not a Markov process unless $\pi(0,\infty)= 0$. \eremark

\btheorem\label{th-HL-Tanaka1b} {\rm(Tanaka's formulae)~} The height process has a local time $\{L_t(a): a, t\ge 0\}$ with stochastic integral representations
 \begin{align}\label{L(a)=bet(H-a)^+2}
L_t(a)= \beta(H_t-a)^+ - \int_0^t 1_{\{H_s> a\}}\d \xi_s + U_t(a)
 \end{align}
and
 \begin{align}\label{L(a)=bet(H-a)^-2}
L_t(a)= \beta(H_t-a)^- + \int_0^t 1_{\{H_s\le a\}}\d \xi_s - I_0(t) - V_t(a),
 \end{align}
where
 \begin{align}\label{U_t(a)=def2}
U_t(a)= \int_0^t\int_0^\infty 1_{\{H_s> a\}} \big(z+I_s(t)-\xi_s\big)^+ N(\d s,\d z)
 \end{align}
and
 \begin{align}\label{V_t(a)=def2}
V_t(a)= \int_0^t\int_0^\infty 1_{\{H_s\le a\}} \big(z+I_s(t)-\xi_s\big)^+ N(\d s,\d z).
 \end{align}
Moreover, there is a version of the local time such that $(a,t)\mapsto L_t(a)$ is jointly measurable and $t\mapsto L_t(a)$ is continuous and increasing. \etheorem

\bproof By \cite[Proposition~3.17]{LPW22}, the local time $L_t(a)$ is given by \eqref{L(a)=bet(H-a)^+2} and \eqref{U_t(a)=def2}. By taking the difference of \eqref{Hproc-intN} and \eqref{L(a)=bet(H-a)^+2}, we obtain \eqref{L(a)=bet(H-a)^-2} and \eqref{V_t(a)=def2}. By \cite[Proposition~1.3.3]{DLG02}, there is a version of the local time such that $(a,t)\mapsto L_t(a)$ is jointly measurable and $t\mapsto L_t(a)$ is continuous and increasing. \eproof

Since the height process possesses a local time, we have $1_{\{H_s=0\}}= 0$ for a.e.\ $s\ge 0$. By taking $a=0$ in \eqref{L(a)=bet(H-a)^-2} we conclude that
 \begin{align}\label{L_t(0)=-I_0(t)}
I_0(t)= -L_t(0), \quad t\ge 0.
 \end{align}
(The above relation is actually true for all $\beta\ge 0$; see \cite[Lemma~1.3.2]{DLG02}.) Similarly, writing \eqref{N=sum_I} for $N(\d s,\d z)$, we have
 \begin{align*}
I_{t_i}(t)-\xi_{t_i}
 =
L_{t_i}(H_{t_i})-L_t(H_{t_i}), \quad t\ge t_i, i\in I.
 \end{align*}
Those calculations lead to some useful stochastic equations given in the following:

\btheorem\label{th-Hproc-inteq} The height process and its local time satisfy the stochastic integral equation
 \begin{align*}
\beta H_t= \xi_t + L_t(0) - \int_0^t \int_0^\infty \big(z + L_s(H_s) - L_t(H_s)\big)^+ N(\d s,\d z), \quad t\ge 0.
 \end{align*}
\etheorem

\bcorollary\label{th-HL-Tanaka2a} {\rm(Tanaka's formulae)~} The height process and its local time satisfy the stochastic integral equations
 \begin{align}\label{L(a)=bet(H-a)^+1}
L_t(a)= \beta(H_t-a)^+ - \int_0^t 1_{\{H_s> a\}}\d \xi_s + U_t(a)
 \end{align}
and
 \begin{align}\label{L(a)=bet(H-a)^-1}
L_t(a)= \beta(H_t-a)^- + \int_0^t 1_{\{H_s\le a\}}\d \xi_s + L_t(0) - V_t(a),
 \end{align}
where
 \begin{align}\label{U_t(a)=def1}
U_t(a)= \int_0^t\int_0^\infty 1_{\{H_s> a\}} \big(z + L_s(H_s) - L_t(H_s)\big)^+ N(\d s,\d z)
 \end{align}
and
 \begin{align}\label{V_t(a)=def1}
V_t(a)= \int_0^t\int_0^\infty 1_{\{H_s\le a\}} \big(z + L_s(H_s) - L_t(H_s)\big)^+ N(\d s,\d z).
 \end{align}
\ecorollary

%\bremark\label{th-casebeta=0} {\blue It does not seem easy to extend the results of this section to the case where $\beta= 0$. For instance, it is not obvious how to get a meaningful limit from the representation \eqref{Hproc-intN} by letting $\beta\to 0$.} \eremark

\section{The spatial $\sigma$-algebras}

\setcounter{equation}{0}

In this section, we introduce a pair of $\sigma$-algebras indexed by the spatial parameter and investigate the measurabilities of some random variables relative to those $\sigma$-algebras. Let $M(\mbb{R}_+)$ denote the space of finite Borel measures on $\mbb{R}_+$ equipped with the topology of weak convergence. Following Duquesne and Le~Gall \cite[p.10 and p.25]{DLG02}, we define the \textit{exploration process} $(\rho_t: t\ge 0)$ by
 \begin{align}\label{exp-proc}
\rho_t(\d x)= \beta 1_{[0,H_t]}(x)\mfr{m}(\d x) + \sum_{{0< t_i\le t}} \itDelta I_{t_i}(t)\delta_{H_s}(\d x),
 \end{align}
where $\mfr{m}(\d x)$ denotes the Lebesgue measure and $\itDelta I_{t_i}(t)$ is given by \eqref{DelI_{t_i}(t)}. The total mass of the exploration process is determined by
 \begin{align}\label{exp-proc-mass}
\<\rho_t,1\>= I_t(t) - I_0(t)= \xi_t + L_t(0), \quad t\ge 0;
 \end{align}
see \cite[p.25 and p.34]{DLG02}.

Let $\supp(\mu)$ denote the topological support of a nonzero measure $\mu\in M(\mbb{R}_+)$. For $\mu\in M(\mbb{R}_+)$ let $H(\mu)= \sup\supp(\mu)$ if $\mu\neq 0$ and let $H(\mu)= 0$ if $\mu= 0$. Then $\mu\mapsto H(\mu)$ is a measurable function on $M(\mbb{R}_+)$. In view of \eqref{exp-proc}, we have
 \begin{align}\label{height-exp-proc}
\mbf{P}\big\{H_t= H(\rho_t)\big\}= 1, \quad t\ge 0.
 \end{align}
Then $\{H(\rho_t): t\ge 0\}$ is a modification of the height process $(H_t: t\ge 0)$; see also Duquesne and Le~Gall \cite[Lemma~1.2.2]{DLG02}.

For a measure $\mu\in M(\mbb{R}_+)$ and a function $f\in C(\mbb{R}_+)$, we write $\<\mu,f\>= \int f\d\mu$. For any $\mu\in M(\mbb{R}_+)$ with compact support, its \textit{concatenation} with $\nu\in M(\mbb{R}_+)$ is the measure $[\mu,\nu]\in M(\mbb{R}_+)$ defined by
 \begin{align*}
\<[\mu,\nu],f\>
 =
\<\mu,f\> + \<\nu,f(H(\mu)+\cdot)\>, \quad f\in C(\mbb{R}_+).
 \end{align*}
For any $a\ge 0$ let $k_a$ be the operator on $M(\mbb{R}_+)$ determined by
 \begin{align*}
k_a\mu([0,x])= \mu([0,x])\land (\<\mu,1\>-a)^+, \quad x\ge 0, \mu\in M(\mbb{R}_+).
 \end{align*}

It was proved by Duquesne and Le~Gall \cite[Proposition~1.2.3]{DLG02} that the exploration process $(\rho_t: t\ge 0)$ defined by \eqref{exp-proc} is a c\`adl\`ag strong Markov process and it is also c\`adl\`ag in the total variation distance. Furthermore, by \cite[formula (20)]{DLG02} for any $t\ge 0$ and any $(\mcr{G}_t)$-stopping time $T$ we have
 \begin{align}\label{rho_{T+t}= [k(t)rho_T,.]}
\rho_{T+t}= \big[k_{-I_0^{(T)}(t)}\rho_T,\rho_t^{(T)}\big],
 \end{align}
where $\xi_0^{(T)}(t)$ and $\rho_t^{(T)}$ denote the analogues of $I_0(t)$ and $\rho_t$ when $\xi_t$ is replaced by $\xi_{T+t} - \xi_T$. By \cite[formula (19)]{DLG02}, the exploration process has the same set of discontinuities as $(\xi_t: t\ge 0)$ and
 \begin{align}\label{rho_t=rho_{t-}+..}
\rho_{t_i}= \rho_{t_i-} + \itDelta \rho_{t_i}
 =
\rho_{t_i-} + \itDelta \xi_{t_i}\delta_{H_{t_i}}, \quad i\in I.
 \end{align}
From \eqref{Hproc-sum0}, \eqref{L_t(0)=-I_0(t)} and \eqref{exp-proc} it is clear that
 \begin{align}\label{<rho_t,1>=xi_t+L_t(0)}
\<\rho_t,1\>= \beta H_t + \sum_{{0< t_i\le t}} \itDelta I_{t_i}(t)= \xi_t + L_t(0).
 \end{align}

For any $a\ge 0$ let $\{\eta_a^+(t): t\ge 0)\}$ and $\{\eta_a^-(t): t\ge 0)\}$ be the positive continuous increasing $(\mcr{G}_t)$-adapted processes defined by
 \begin{align*}
\eta_a^+(t)= \int_0^t 1_{\{H_s> a\}}\d s,
 \quad
\eta_a^-(t)= \int_0^t 1_{\{H_s\le a\}}\d s.
 \end{align*}
Clearly, we have $\eta_a^+(t)\to \infty$ for $a\ge 0$ and $\eta_a^-(t)\to \infty$ for $a> 0$ as $t\to \infty$. For $a\ge 0$ and $s\ge 0$ define the $(\mcr{G}_t)$-stopping times
 \begin{align*}
\tau_a^+(s)= \inf\{t\ge 0: \eta_a^+(t)> s\},
 \quad
\tau_a^-(s)= \inf\{t\ge 0: \eta_a^-(t)> s\}.
 \end{align*}
Then $s\mapsto \tau_a^+(s)$ is a left-continuous strictly increasing process for $a\ge 0$ and $s\mapsto \tau_a^-(s)$ is a left-continuous strictly increasing process for $a> 0$. Note that $\tau_0^-(s)= \infty$ for $s\ge 0$. For any $a\ge 0$ define the $M(\mbb{R}_+)$-valued process
 \begin{align}\label{rho^a=def}
\<\rho_t^a,f\>= \int_a^\infty f(x-a)\rho_{\tau_a^+(t)}(\d x), \quad t\ge 0.
 \end{align}
Let $\mcr{F}_0$ be the $\sigma$-algebra generated by the collection of $\mbf{P}$-null sets. For $a> 0$ let $\mcr{F}_a$ be the $\sigma$-algebra generated by $\mcr{F}_0$ and the process $(\rho_{\tau_a^-(t)}: t\ge 0)$. For $a\ge 0$ let $\mcr{F}^a$ be the $\sigma$-algebras generated by $\mcr{F}_0$ and the process $(\rho_{\tau_a^+(t)}: t\ge 0)$. The $\sigma$-algebras $\mcr{F}^a$ and $\mcr{F}_a$ contain the information of the random noises of $(\xi_t: t\ge 0)$ during the times when the height process is above and below the level $a\ge 0$, respectively. For $a\ge 0$ let
 \begin{align*}
\mcr{F}_{a+}= \bigcap_{c> a}\mcr{F}_c.
 \end{align*}

\bproposition\label{explor-indincr} For any $a\ge 0$ the $\sigma$-algebras $\mcr{F}_{a+}$ and $\mcr{F}^a$ are independent. \eproposition

\bproof By \eqref{rho^a=def} we see that $\mcr{F}^a$ is also generated by $\{\rho_t^a: t\ge 0\}$, so it is independent of $\mcr{F}_a$ by \cite[Proposition~1.3.1]{DLG02}. Then $\mcr{F}_{a+}$ and $\mcr{F}^c$ are independent for any $c> a$. A monotone class argument shows that $\mcr{F}_{a+}$ and $\mcr{F}^a= \sigma(\cup_{c> a}\mcr{F}_c)$ are independent. \eproof

We next introduce some time changed processes and discuss their measurabilities relative to the $\sigma$-algebras $\mcr{F}^a$ and $\mcr{F}_a$ for $a\ge 0$. Let $\mcr{G}_{a,t}^{\pm}= \mcr{G}_{\tau_a^{\pm}(t)}$ for any $a\ge 0$ and $t\ge 0$. Write
 \begin{align*}
B_{a,t}^+= \int_0^{\tau_a^+(t)} 1_{\{a< H_s\}}\d B_s,
 ~~
B_{a,t}^-= \int_0^{\tau_a^-(t)} 1_{\{H_s\le a\}} \d B_s.
 \end{align*}
Then $(B_{a,t}^+: t\ge 0)$ is a standard $(\mcr{G}_{a,t}^+)$-Brownian motion for $a\ge 0$ and $(B_{a,t}^-: t\ge 0)$ is a standard $(\mcr{G}_{a,t}^-)$-Brownian motion for $a> 0$. Note that $B_{0,t}^-= 0$ a.s.\ for each $t\ge 0$. Clearly, for any $a\ge 0$ the processes $(B_{a,t}^+: t\ge 0)$ and $(B_{a,t}^-: t\ge 0)$ are independent.

For $a\ge 0$ we can define a $(\mcr{G}_{a,t}^+)$-Poisson random measure $N_a^+(\d s,\d z)$ on $(0,\infty)^2$ with intensity $\d s\pi(\d z)$ by
 \begin{align}\label{N_a^+=def1}
N_a^+((0,t]\times B)= \int_0^{\tau_a^+(t)}\int_B 1_{\{a< H_s\}} N(\d s,\d z).
 \end{align}
Similarly, for $a> 0$ we can define a $(\mcr{G}_{a,t}^-)$-Poisson random measure $N_a^-(\d s,\d z)$ on $(0,\infty)^2$ with intensity $\d s\pi(\d z)$ by
 \begin{align}\label{N_a^-=def1}
N_a^-((0,t]\times B)= \int_0^{\tau_a^-(t)}\int_B 1_{\{H_s\le a\}} N(\d s,\d z).
 \end{align}
For completeness of the notation, we denote by $N_0^-(\d s,\d z)$ the trivial measure on $(0,\infty)^2$, that is $N_0^-((0,\infty)^2)= 0$. Then for any $a\ge 0$ the random measures $N_a^+(\d s,\d z)$ and $N_a^-(\d s,\d z)$ are independent.

For any $a\ge 0$ we define the $(\mcr{G}_{a,t}^+)$-L\'evy processes $(\xi_{a,t}^+: t\ge 0)$ and the $(\mcr{G}_{a,t}^-)$-L\'evy processes $(\xi_{a,t}^-: t\ge 0)$ by
 \begin{align}\label{xi_t^{a+-}=def2}
\xi_{a,t}^+= \int_0^{\tau_a^+(t)} 1_{\{a< H_s\}}\d \xi_s,
 \quad
\xi_{a,t}^-= \int_0^{\tau_a^-(t)} 1_{\{H_s\le a\}} \d \xi_s.
 \end{align}
It is easy to see that
 \begin{align}\label{xi_t^{a+-}=def1}
\xi_{a,t}^{\pm}= - \alpha t + \sqrt{2\beta}B_{a,t}^{\pm} + \int_0^t\int_0^\infty z \tilde{N}_a^{\pm}(\d s,\d z).
 \end{align}
Then $(\xi_{a,t}^+: t\ge 0)$ and $(\xi_{a,t}^-: t\ge 0)$ are independent.

Let $H_{a,t}^{\pm}= H_{\tau_a^{\pm}(t)}$ for $a\ge 0$ and $t\ge 0$. Then $(H_{a,t}^+: t\ge 0)$ is a continuous $(\mcr{G}_{a,t}^+)$-adapted process and $(H_{a,t}^-: t\ge 0)$ is a continuous $(\mcr{G}_{a,t}^-)$-adapted process. It is simple to see that $(H_{a,t}^+: t\ge 0)$ has local time $(L_{a,t}^+(h): t\ge 0)$ at level $h\ge a$ given by
 \begin{align}\label{L_{a,t}^+(a)=..}
L_{a,t}^+(h):= \int_0^{\tau_a^+(t)} 1_{\{a< H_s\}}\d L_s(h).
 \end{align}
Similarly, the process $(H_{a,t}^-: t\ge 0)$ has local time $(L_{a,t}^-(h): t\ge 0)$ at level $h\le a$ given by
 \begin{align}\label{L_{a,t}^-(a)=..}
L_{a,t}^-(h):= \int_0^{\tau_a^-(t)} 1_{\{H_s\le a\}}\d L_s(h).
 \end{align}

\bproposition\label{thL_{tau_a(t)}=L_{a,t}} {\rm(i)} For any $h\ge a\ge 0$ we have
 \begin{align*}
\mbf{P}\big(L_{\tau_a^+(t)}(h)= L_{a,t}^+(h) ~ \mbox{and}~  L_t(h)= L_{a,\eta_a^+(t)}^+(h) ~ \mbox{for all}~ t\ge 0\big)= 1.
 \end{align*}

{\rm(ii)} For any $a\ge h\ge 0$ we have
 \begin{align*}
\mbf{P}\big(L_{\tau_a^-(t)}(h)= L_{a,t}^-(h) ~ \mbox{and}~  L_t(h)= L_{a,\eta_a^-(t)}^-(h) ~ \mbox{for all}~ t\ge 0\big)= 1.
 \end{align*}
\eproposition

\bproof We only give the proof of (ii) since (i) follows by similar arguments. In this proof, we may assume $a> 0$, for otherwise the result holds trivially. By \cite[Proposition~1.3.3]{DLG02} we have, for every $v\ge 0$,
 \begin{align*}%\label{limsupE(1_{h+vep})=0}
\lim_{\vepsilon\to 0+}\sup_{h\ge 0}\mbf{E}\Big[\sup_{0\le t\le v} \Big|\frac{1}{\vepsilon} \int_0^t 1_{\{h< H_s\le h+\vepsilon\}}\d s - L_t(h)\Big|\Big]= 0
 \end{align*}
and
 \begin{align*}%\label{limsupE(1_{h-vep})=0}
\lim_{\vepsilon\to 0+}\sup_{h\ge \vepsilon}\mbf{E}\Big[\sup_{0\le t\le v} \Big|\frac{1}{\vepsilon} \int_0^t 1_{\{h-\vepsilon< H_s\le h\}}\d s - L_t(h)\Big|\Big]= 0.
 \end{align*}
It follows that for any $h\ge 0$ there is a sequence $\{\vepsilon_k\}\subset (0,1]$ decreasing to zero such that the following property a.s.\ holds: for all $t\ge 0$,
 \begin{align}\label{L_t(h)=limfrac{1}{vep_k}(+)0}
L_t(h)= \lim_{k\to \infty} \frac{1}{\vepsilon_k} \int_0^t 1_{\{h< H_s\le h+\vepsilon_k\}}\d s.
 \end{align}
Similarly, for any $h> 0$ there is a sequence $\{\vepsilon_k\}\subset (0,1]$ decreasing to zero such that the following property a.s.\ holds: for all $t\ge 0$,
 \begin{align}\label{L_t(h)=limfrac{1}{vep_k}(-)0}
L_t(h)= \lim_{k\to \infty} \frac{1}{\vepsilon_k} \int_0^t 1_{\{h-\vepsilon_k< H_s\le h\}}\d s.
 \end{align}
By \eqref{L_t(h)=limfrac{1}{vep_k}(-)0} for $a\ge h> 0$ the following property a.s.\ holds: for all $t\ge 0$,
 \begin{align*}%\label{L(h)=limfrac(1)}
L_{\tau_a^-(t)}(h)\aar= \lim_{k\to \infty} \frac{1}{\vepsilon_k} \int_0^{\tau_a^-(t)} 1_{\{h-\vepsilon_k< H_s\le h\}}\d s \cr
 \aar=
\lim_{k\to \infty} \frac{1}{\vepsilon_k} \int_0^t 1_{\{h-\vepsilon_k< H_{a,s}^-\le h\}}\d s= L_{a,t}^-(h)
 \end{align*}
and
 \begin{align*}%\label{L(h)=limfrac(2)}
L_t(h)\aar= \lim_{k\to \infty} \frac{1}{\vepsilon_k} \int_0^t 1_{\{h< H_s\le h+\vepsilon_k\}}\d s \cr
 \aar=
\lim_{k\to \infty} \frac{1}{\vepsilon_k} \int_0^{\eta_a^-(t)} 1_{\{h-\vepsilon_k< H_{a,s}^-\le h\}}\d s= L_{a,\eta_a^-(t)}^-(h).
 \end{align*}
Taking $h=0$ in \eqref{L_t(h)=limfrac{1}{vep_k}(+)0} we see that the following property a.s.\ holds: for all $t\ge 0$,
 \begin{align*}%\label{L(0)=limfrac(1)}
L_{\tau_a^-(t)}(0)\aar= \lim_{k\to \infty} \frac{1}{\vepsilon_k} \int_0^{\tau_a^-(t)} 1_{\{0< H_s\le \vepsilon_k\}}\d s \cr
 \aar=
\lim_{k\to \infty} \frac{1}{\vepsilon_k} \int_0^t 1_{\{0< H_{a,s}^-\le \vepsilon_k\}}\d s= L_{a,t}^-(0)
 \end{align*}
and
 \begin{align*}%\label{L(0)=limfrac(2)}
L_t(0)\aar= \lim_{k\to \infty} \frac{1}{\vepsilon_k} \int_0^t 1_{\{0< H_s\le \vepsilon_k\}}\d s \cr
 \aar=
\lim_{k\to \infty} \frac{1}{\vepsilon_k} \int_0^{\eta_a^-(t)} 1_{\{0< H_{a,s}^-\le \vepsilon_k\}}\d s= L_{a,\eta_a^-(t)}^-(0).
 \end{align*}
Then the desired result follows. \eproof

\bproposition\label{th(H^{a+-})-msbl} {\rm(i)} For any $a\ge 0$ the following random variables are $\mcr{F}^a$-measurable:
 \begin{align*}
\xi_{a,t}^+, H_{a,t}^+, L_{a,t}^+(h), N_a^+(A), \quad t\ge 0, h\ge a, A\in \mcr{B}((0,\infty)^2).
 \end{align*}
Furthermore, if $\beta> 0$, then the random variables $B_{a,t}^+$, $t\ge 0$ are $\mcr{F}^a$-measurable.

{\rm(ii)} For any $a\ge 0$ the following random variables are $\mcr{F}_a$-measurable:
 \begin{align*}
\xi_{a,t}^-, H_{a,t}^-, L_{a,t}^-(h), N_a^-(A), \quad t\ge 0, 0\le h\le a, A\in \mcr{B}((0,\infty)^2).
 \end{align*}
Furthermore, if $\beta> 0$, then the random variables $B_{a,t}^-$, $t\ge 0$ are $\mcr{F}_a$-measurable. \eproposition

\bproof We only give the proof of (ii) since (i) follows similarly. By \eqref{height-exp-proc} and the definition of $\mcr{F}_a$ it is clear that $H_{a,t}^-$ is $\mcr{F}_a$-measurable. Then the local time $L_{a,t}^-(h)= L_{\tau_a^-(t)}(h)$ is $\mcr{F}_a$-measurable for $a\ge h\ge 0$. From \eqref{rho_t=rho_{t-}+..} and \eqref{N_a^-=def1} we see that $N_a^-(A)$ is $\mcr{F}_a$-measurable. If $\beta= 0$, then $\xi_{a,t}^-$ is $\mcr{F}_a$-measurable by \eqref{xi_t^{a+-}=def1}. When $\beta> 0$, we get from \eqref{L(a)=bet(H-a)^-1} and Proposition~\ref{thL_{tau_a(t)}=L_{a,t}}-(ii) that
 \begin{align*}
L_{a,t}^-(a)= \beta(H_{a,t}^--a)^- + \xi_{a,t}^- + L_{a,t}^-(0) - V_{a,t}^-(a),
 \end{align*}
where
 \begin{align*}
V_{a,t}^-(a)= \int_0^t\int_0^\infty \big(z + L_{a,s}^-(H_{a,s}^-) - L_{a,t}^-(H_{a,s}^-)\big)^+ N_a^-(\d s,\d z).
 \end{align*}
Then $\xi_{a,t}^-$ is again $\mcr{F}_a$-measurable, and the $\mcr{F}_a$-measurability of $B_{a,t}^-$ follows by \eqref{xi_t^{a+-}=def1}. \eproof

\bremark\label{thmcr{F}_a=mcr{H}_a} The definition of the spatial filtration $(\mcr{F}_a: a\ge 0)$ is slightly different from that of $(\mcr{H}_a: a\ge 0)$ given in \cite[p.31]{DLG02}, where $\mcr{H}_a$ was generated by $\{(\xi_{\tau_a^-(t)}, \rho_{\tau_a^-(t)}): t\ge 0\}$. However, by \eqref{exp-proc-mass} and Proposition~\ref{thL_{tau_a(t)}=L_{a,t}}-(ii) we have
 \begin{align*}
\xi_{\tau_a^-(t)}= \<\rho_{\tau_a^-(t)},1\> - L_{\tau_a^-(t)}(0)
 =
\<\rho_{\tau_a^-(t)},1\> - L_{a,t}^-(0),
 \end{align*}
which is $\mcr{F}_a$-measurable by Proposition~\ref{th(H^{a+-})-msbl}-(ii). Then we actually have $\mcr{F}_a= \mcr{H}_a$ for all $a\ge 0$. \eremark

\section{The local time process}

\setcounter{equation}{0}

From Duquesne and Le~Gall \cite[Theorem~1.4.1]{DLG02} we know that $\{L_{T_x}(a): a\ge 0\}$ is a CB-process relative to the filtration $(\mcr{F}_a: a\ge 0)$ with transition semigroup defined by \eqref{LaplaceQ_t(x,.)}. Then it has a c\`adl\`ag modification, which is a CB-process relative to the right-continuous filtration $(\mcr{F}_{a+}: a\ge 0)$ with the same transition semigroup. Moreover, we have
 \begin{align}\label{E[L_{T_x}(a)]=xe^{}}
\mbf{E}[L_{T_x}(a)]= x\e^{-\alpha a}, \quad a\ge 0.
 \end{align}

We shall give the proof of Theorem~\ref{thL_{T_x}(a)=x+int..} in the cases $\beta> 0$ and $\beta= 0$ separately. In the first case, we obtain \eqref{L_{T_x}(t)=x+int_0^{T_x}} from \eqref{L(a)=bet(H-a)^-2}. In the second case, the proof of \eqref{L_{T_x}(t)=x+int_0^{T_x}} is based on a semimartingale characterizations of the CB-process $\{L_{T_x}(a): a\ge 0\}$ given in \cite{Li20}.

\blemma\label{th1_{sleT_x}=a.s.} The left-continuous $(\mcr{G}_t)$-adapted process $s\mapsto 1_{\{s\le T_x\}}$ is a modification of $s\mapsto 1_{\{L_s(H_s)\le L_{T_x}(H_s)\}}$. \elemma

\bproof For any $v,u> 0$ set $g_*(v,u)= 1_{\{u< L_{T_x}(v)\}}$ and $g^*(v,u)= 1_{\{u\le L_{T_x}(v)\}}$. Since the local time process $t\mapsto L_t(h)$ is increasing, we have
 \begin{align*}
\{L_s(H_s)< L_{T_x}(H_s)\}\subset \{s< T_x\}\subset \{s\le T_x\}\subset \{L_s(H_s)\le L_{T_x}(H_s)\},
 \end{align*}
and hence
 \begin{align*}
\qquad g_*(H_s,L_s(H_s))\le 1_{\{s< T_x\}}
 \le
1_{\{s\le T_x\}}\le g^*(H_s,L_s(H_s)).
 \end{align*}
By the occupation times formula, %(see, e.g., Revuz and Yor \cite[p.232]{ReY05})
 \begin{align*}
\aar\mbf{E}\Big\{\int_0^\infty [g^*(H_s,L_s(H_s)) - g_*(H_s,L_s(H_s))]\d s\Big\} \cr
 \aar\qqquad
= \mbf{E}\Big\{\int_0^\infty\d v\int_0^\infty [g^*(v,L_s(v)) - g_*(v,L_s(v))]\d L_s(v)\Big\} \cr
 \aar\qqquad
= \mbf{E}\Big\{\int_0^\infty\d v\int_0^\infty [g^*(v,u) - g_*(v,u)]\d u\Big\} \cr
 \aar\qqquad
= \mbf{E}\Big\{\int_0^\infty\d v\int_0^\infty 1_{\{u= L_{T_x}(v)\}}\d u\Big\}= 0,
 \end{align*}
which implies
 \begin{align*}
\mbf{E}\Big\{\int_0^\infty [g^*(H_s,L_s(H_s)) - 1_{\{s\le T_x\}}]\d s\Big\}= 0.
 \end{align*}
That proves the desired result. \eproof

\medskip\noindent\textit{Proof of Theorem~\ref{thL_{T_x}(a)=x+int..} (for $\beta> 0$).~} The definition of the hitting time $T_x$ implies that $L_{T_x}(0)= x$ and $H_{T_x}= 0$. Clearly, for any $i\in I$ the process $t\mapsto L_t(0)$ does not increase in the interval $(t_i,r(t_i))$, and hence $T_x\notin (t_i,r(t_i))$. If $r(t_i)< T_x$, then
 \begin{align*}
\big(z_i + I_{t_i}(T_x) - \xi_{t_i}\big)^+
 =
\big(z - L_{T_x}(H_{t_i}) + L_{t_i}(H_{t_i})\big)^+
 =
\big(z_i - \xi_{t_i}\big)^+= 0.
 \end{align*}
Then taking $t= T_x$ in \eqref{L(a)=bet(H-a)^-2} or \eqref{L(a)=bet(H-a)^-1} we obtain \eqref{L_{T_x}(t)=x+int_0^{T_x}}. \qed

\medskip\noindent\textit{Proof of Theorem~\ref{thL_{T_x}(a)=x+int..} (for $\beta= 0$).~} To avoid triviality, we assume $\pi(0,\infty)\neq 0$. Suppose that the Poisson random measure $N(\d s,\d z)$ is given by \eqref{N=sum_I}.

\textit{Step~1.} By Li \cite[Theorem~7.2]{Li20}, the CB-process $\{L_{T_x}(a): a\ge 0\}$ is a semimartingale with decomposition
 \begin{align}\label{L_{T_x}(a)=x-alp+(2)}
L_{T_x}(a)= x - \alpha\int_0^a L_{T_x}(v) \d v + \int_0^a\int_0^\infty z \tilde{M}_x(\d v,\d z),
 \end{align}
where $\tilde{M}_x(\d s,\d z)$ is an optional random measure with predictable compensator
 \begin{align}\label{hat{M}_x(ds,dz)=..}
\hat{M}_x(\d s,\d z)= 1_{\{z> \vepsilon\}}L_{T_x}(s-)\d s\pi(\d z)
 =
1_{\{z> \vepsilon\}}L_{T_x}(s)\d s\pi(\d z).
 \end{align}
By the occupation time formula,
 \begin{align*}
\int_0^a L_{T_x}(v) \d v= \int_0^{T_x} 1_{\{H_s\le a\}}\d s.
 \end{align*}

\textit{Step~2.} Take $\vepsilon> 0$ such that $\pi(\vepsilon,\infty)> 0$. Let $s_1= \inf\{t\ge 0: \xi_t - \xi_{t-}> \vepsilon\}$ and $r(s_1)= \inf\{t\ge s_1: \xi_t= \xi_{s_1-}\}$. For $i\ge 2$ define inductively $s_i= \inf\{t\ge r(s_{i-1}): \xi_t - \xi_{t-}> \vepsilon\}$ and $r(s_i)= \inf\{t\ge s_i: \xi_t= \xi_{s_i-}\}$. Let $D_\vepsilon= \cup_{k=1}^\infty (s_i,r(s_i)]$ and let
 \begin{align}
\eta_\vepsilon(t)= \int_0^t 1_{D_\vepsilon^c}(s)\d s, \quad t\ge 0,
 \end{align}
where $D_\vepsilon^c= [0,\infty)\setminus D_\vepsilon$. Let $\gamma_\vepsilon(t)= \inf\{s\ge 0: \eta_\vepsilon(s)> t\}$ and
 \begin{align*}
\xi_t^{(0)}= \xi_{\gamma_\vepsilon(t)}= \int_0^{\gamma_\vepsilon(t)} 1_{D_\vepsilon^c}(s)\d \xi_s, \quad t\ge 0.
 \end{align*}
Then $\{\xi_t^{(0)}: t\ge 0\}$ is a L\'evy process with Laplace exponent $\psi_\vepsilon$ given by
 \begin{align}\label{psi_vep(lam)=def}
\psi_\vepsilon(\lambda)= \Big[\alpha + \int_\vepsilon^\infty z\pi(\d z)\Big]\lambda + \int_0^\vepsilon (\e^{-\lambda z}-1+\lambda z) \pi(\d z).
 \end{align}
Let $T_x^{(0)}= \inf\{t\ge 0: \xi_t^{(0)}= -x\}$. Observe that $\gamma_\vepsilon(T_x^{(0)})= T_x$. Let $\{H_t^{(0)}: t\ge 0\}$ be the height process of $\{\xi_t^{(0)}: t\ge 0\}$. Then the corresponding local time process $\{L_{T_x^{(0)}}(a,H^{(0)}): a\ge 0\}$ is a CB-process with branching mechanism $\psi_\vepsilon$ given by \eqref{psi_vep(lam)=def}.

\textit{Step~3.} Let $\xi_t^{(i)}= \xi_{s_i+t} - \xi_{s_i}$. Then $\{\xi_t^{(i)}: t\ge 0\}$ is also a L\'evy process with Laplace exponent $\psi$. Let $\{H^{(i)}(t): t\ge 0\}$ be the corresponding height process and $\{\rho_t^{(i)}: t\ge 0\}$ be the corresponding exploration process. Let $\{L_t^{(i)}(a): a,t\ge 0\}$ be the local time of $\{H_t^{(i)}: t\ge 0\}$. In view of \eqref{rho_t=rho_{t-}+..}, we have
 \begin{align*}
\rho_{s_i}= \rho_{s_i-} + \itDelta \rho_{s_i}
 =
\rho_{s_i-} + \itDelta \xi_{s_i}\delta_{H_{s_i}}.
 \end{align*}
From \eqref{rho_{T+t}= [k(t)rho_T,.]} it follows that
 \begin{align}
\rho_{s_i+t}= \big[k_{-I_0^{(i)}(t)}\rho_{s_i},\rho_t^{(i)}\big], \quad t\ge 0.
 \end{align}
Recall that $H_t= \sup\supp(\rho_t)$. Then the above relation implies that
 \begin{align*}
H_{s_i+t}= H_{s_i} + H^{(i)}(t), \quad t\ge 0.
 \end{align*}
By \cite[Lemma~1.3.2]{DLG02} we see that
 \begin{align}
L_t^{(i)}(0)= I_0^{(i)}(t):= \inf_{0\le s\le t} \xi_s^{(i)}.
 \end{align}
It is clear that $r(s_i)-s_i= \inf\{t\ge 0: \xi_t^{(i)}= -z_i\}$. Then $\{L_{r(s_i)-s_i}^{(i)} (a): a\ge 0\}$ is a CB-process with branching mechanism $\psi$ and $L_{r(s_i)-s_i}^{(i)}(0)= z_i$. Observe that
 \begin{align}\label{L_{T_x}(a,H)=..}
L_{T_x}(a)= L_{T_x^{(0)}}^{(0)}(a) + \sum_{s_i\le T_x} 1_{\{H_{s_i}\le a\}} L_{r(s_i)-s_i}^{(i)} (a-H_{s_i}).
 \end{align}

\textit{Step~4.} Clearly, each term on the r.h.s.\ of \eqref{L_{T_x}(a,H)=..} can only make positive jumps as $a\ge 0$ increases. Therefore, on the event $\{s_i\le T_x\}$ we have
 \begin{align*}
\Delta L_{T_x}(H_{s_i}):= L_{T_x}(H_{s_i})-L_{T_x}(H_{s_i}-)
 \ge
L_{r(s_i)-s_i}^{(i)}(0)= z_i.
 \end{align*}
Let $\vepsilon> 0$ be given as in the second step. It follows that
 \begin{align}\label{int1_{}zN(dsdz)le..}
\int_0^a\int_\vepsilon^\infty z M_x(\d v,\d z)
 \aar=
\sum_{v\le a} 1_{\{\Delta L_{T_x}(v)> \vepsilon\}}\Delta L_{T_x}(v) \cr
 \aar\ge
\sum_{s_i\le T_x} 1_{\{H_{s_i}\le a, z_i> \vepsilon\}} z_i \cr
 \aar=
\int_0^{T_x}\int_\vepsilon^\infty 1_{\{H_s\le a\}} z N(\d s,\d z).
 \end{align}
In view of \eqref{E[L_{T_x}(a)]=xe^{}} and \eqref{hat{M}_x(ds,dz)=..}, we have
 \begin{align*}
\mbf{E}\Big[\int_0^a\int_\vepsilon^\infty z M_x(\d v,\d z)\Big]
 \aar=
\mbf{E}\Big[\int_0^a L_{T_x}(v) \d v\int_\vepsilon^\infty z \pi(\d z)\Big] \cr
 \aar=
x\int_0^a\e^{-\alpha v} \d v\int_\vepsilon^\infty z \pi(\d z).
 \end{align*}
By Lemma~\ref{th1_{sleT_x}=a.s.} and the occupation time formula,
 \begin{align*}
\int_0^\infty 1_{\{s\le T_x, H_s\le a\}}\d s
 \aar=
\int_0^\infty 1_{\{L_s(H_s)\le L_{T_x}(H_s), H_s\le a\}} \d s \cr
 \aar=
\int_0^a\d v\int_0^\infty 1_{\{L_s(v)\le L_{T_x}(v)\}}\d L_s(v) \cr
 \aar=
\int_0^a\d v\int_0^\infty 1_{\{u\le L_{T_x}(v)\}}\d u \cr
 \aar=
\int_0^a L_{T_x}(v) \d v.
 \end{align*}
Consequently, we can use \eqref{E[L_{T_x}(a)]=xe^{}} to see that
 \begin{align*}
\mbf{E}\Big[\int_0^{T_x}\int_\vepsilon^\infty 1_{\{H_s\le a\}} z N(\d s,\d z)\Big]
 \aar=
\mbf{E}\Big[\int_0^\infty 1_{\{s\le T_x, H_s\le a\}}\d s\int_\vepsilon^\infty z \pi(\d z)\Big] \cr
 \aar=
\mbf{E}\Big[\int_0^a L_{T_x}(v) \d v\int_\vepsilon^\infty z \pi(\d z)\Big] \cr
 \aar=
x\int_0^a\e^{-\alpha v} \d v\int_\vepsilon^\infty z \pi(\d z).
 \end{align*}
Then \eqref{int1_{}zN(dsdz)le..} implies that
 \begin{align*}
\int_0^a\int_\vepsilon^\infty z M_x(\d v,\d z)
 =
\int_0^{T_x}\int_\vepsilon^\infty 1_{\{H_s\le a\}} z N(\d s,\d z).
 \end{align*}
Since $\vepsilon> 0$ could be arbitrarily small, we get \eqref{L_{T_x}(t)=x+int_0^{T_x}} from \eqref{L_{T_x}(a)=x-alp+(2)}. \qed

%\bremark\label{proofforbeta=0} {\blue It is interesting to give a proof of Theorem~\ref{thL_{T_x}(a)=x+int..} in the case $\beta= 0$ without using the fact that $\{L_{T_x}(a): a\ge 0\}$ is a CB-process, which would lead to a new proof of the Ray--Knight theorem.} \eremark

\section{The time-space noises}

\setcounter{equation}{0}

In this section, we discuss the construction and characterization of the time-space L\'evy noise provided in Proposition~\ref{th-GP-noises}. Recall that $\mcr{L}_2$ and $\mcr{L}_3(\pi)$ are the function spaces defined by \eqref{iintf<infty} and \eqref{iiintg<infty}, respectively. For $f\in \mcr{L}_2$ and $g\in \mcr{L}_3(\pi)$, define the two-parameter processes $\{W_t(a,f): a,t\ge 0\}$ and $\{M_t(a,g): a,t\ge 0\}$ by
 \begin{align}\label{W_t(a,f)=def}
W_t(a,f)= \int_0^t 1_{\{H_s\le a\}} f(H_s,L_s(H_s))\d B_s
 \end{align}
and
 \begin{align}\label{tilM_t(a,f)=def}
\tilde{M}_t(a,f)= \int_0^t\int_0^\infty 1_{\{H_s\le a\}} g(H_s,z,L_s(H_s)) \tilde{N}(\d s,\d z).
 \end{align}
Here we should note that the process $t\mapsto L_t(H_t)$ is $(\mcr{G}_t)$-predictable since both $t\mapsto H_t$ and $t\mapsto L_t(a)$ are continuous and $(\mcr{G}_t)$-adapted processes.

For any Lebesgue integrable function $f$ on $[0,\infty)^2$, one can use the occupation times formula to see that
 \begin{align}\label{intf(H_s,L_s(H_s))=}
\int_0^\infty f(H_s,L_s(H_s))\d s\aar= \int_0^\infty\d v\int_0^\infty f(v,L_s(v))\d L_s(v) \cr
 \aar=
\int_0^\infty \d v\int_0^\infty f(v,u)\d u.
 \end{align}

\bproposition\label{thW(a,f)=def} For any $a\ge 0$ and $f\in \mcr{L}_2$ the limit $W(a,f):= \lim_{t\to \infty} W_t(a,f)$ exists almost surely and in the $L^2$ sense, and it is given by
 \begin{align}\label{W(a,f)=def}
W(a,f)= \int_0^\infty 1_{\{H_s\le a\}} f(H_s,L_s(H_s))\d B_s,
 \end{align}
which is a centered Gaussian random variable with variance
 \begin{align}\label{<W(a,f)>=iintf^2}
\<W(a,f)\>= \int_0^a \d s\int_0^\infty f(s,u)^2\d u.
 \end{align}
\eproposition

\bproof The process $\{W_t(a,f): t\ge 0\}$ is a continuous square integrable $(\mcr{G}_t)$-martingale with quadratic variation process given by
 \begin{align*}
\<W(a,f)\>_t= \int_0^t 1_{\{H_s\le a\}}f(H_s,L_s(H_s))^2\d s.
 \end{align*}
Clearly, we have the increasing limit
 \begin{align*}
\lim_{t\to \infty}\<W(a,f)\>_t= \<W(a,f)\>:= \int_0^\infty 1_{\{H_s\le a\}}f(H_s,L_s(H_s))^2\d s.
 \end{align*}
In view of \eqref{intf(H_s,L_s(H_s))=}, we can write $\<W(a,f)\>$ as in \eqref{<W(a,f)>=iintf^2}. Then the limit $W(a,f):= \lim_{t\to \infty} W_t(a,f)$ exists almost surely and in the $L^2$ sense, and it can be written as the stochastic integral \eqref{W(a,f)=def}. By It\^o's formula, one can see that
 \begin{align}\label{e^{iW_t(a,f)+..}}
\e^{\im W_t(a,f)+\frac{1}{2}\<W(a,f)\>_t}= 1 + \im\int_0^t \e^{\im W_s(a,f)+\frac{1}{2}\<W(a,f)\>_s} 1_{\{H_s\le a\}} f(H_s,L_s(H_s)) \d B_s,
 \end{align}
which is a continuous square integrable $(\mcr{G}_t)$-martingale. It follows that
 \begin{align*}
\mbf{E}\big[\e^{\im W_t(a,f)+\frac{1}{2}\<W(a,f)\>_t}\big]= 1.
 \end{align*}
By letting $t\to \infty$ in the above equality and using dominated convergence we obtain, for any $a\ge 0$ and $f\in \mcr{L}_2$,
 \begin{align*}
\mbf{E}\big[\e^{\im W(a,f)+\frac{1}{2}\<W(a,f)\>}\big]= 1.
 \end{align*}
Then $W(a,f)$ is a centered Gaussian random variable with variance $\<W(a,f)\>$. \eproof

\bproposition\label{thtilM(a,g)=def} For any $a\ge 0$ and $g\in \mcr{L}_3(\pi)$ the limit $\tilde{M}(a,g):= \lim_{t\to \infty} \tilde{M}_t(a,g)$ exists almost surely and in the $L^2$ sense, and it is given by
 \begin{align}\label{tilM(a,g)=def}
\tilde{M}(a,g)= \int_0^\infty\int_0^\infty 1_{\{H_s\le a\}} g(H_s,z,L_s(H_s)) \tilde{N}(\d s,\d z).
 \end{align}
Moreover, we have
 \begin{align}\label{E[e^{iM(a,g)}]=}
\mbf{E}\big[\e^{\im \tilde{M}(a,g)}\big]
 =
\e^{A(a,g)}, \quad a\ge 0, g\in \mcr{L}_3(\pi),
 \end{align}
where
 \begin{align*}
A(a,g)= \int_0^a\d s\int_0^\infty \pi(\d z) \int_0^\infty \big[\e^{\im g(s,z,u)} - 1 - \im g(s,z,u)\big]\d u.
 \end{align*}
\eproposition

\bproof By \eqref{tilM_t(a,f)=def} one can see that $\{\tilde{M}_t(a,g): t\ge 0\}$ is a c\`adl\`ag and square integrable $(\mcr{G}_t)$-martingale with
 \begin{align*}
\mbf{E}\big[\tilde{M}_t(a,g)^2\big]= \int_0^t\d s\int_0^\infty 1_{\{H_s\le a\}} g(H_s,z,L_s(H_s))^2 \pi(\d z).
 \end{align*}
In view of \eqref{intf(H_s,L_s(H_s))=}, we have
 \begin{align*}
\mbf{E}\big[\tilde{M}_t(a,g)^2\big]\aar\le \int_0^\infty 1_{\{H_s\le a\}} \d s\int_0^\infty g(H_s,z,L_s(H_s))^2 \pi(\d z) \cr
 \aar=
\int_0^a\d s\int_0^\infty \pi(\d z) \int_0^\infty g(s,z,u)^2\d u< \infty.
 \end{align*}
By a martingale convergence theorem, the limit $\tilde{M}(a,g):= \lim_{t\to \infty} \tilde{M}_t(a,g)$ exists almost surely and in the $L^2$ sense. Clearly, the limit can be written as the stochastic integral \eqref{tilM(a,g)=def}. For $t\ge 0$ let
 \begin{align}\label{A_t(a,g)=def}
A_t(a,g)= \int_0^t\d s\int_0^\infty 1_{\{H_s\le a\}} K_g(H_s,z,L_s(H_s)) \pi(\d z).
 \end{align}
By It\^o's formula, we have{\small
 \begin{align}\label{e^{itilM_t(a,f)-..}}
\e^{\im \tilde{M}_t(a,g)-A_t(a,g)}\aar= 1 + \im\int_0^t\int_0^\infty \e^{\im\tilde{M}_{s-}(a,g)-A_s(a,g)} 1_{\{H_s\le a\}} g(H_s,z,L_s(H_s)) \tilde{N}(\d s,\d z) \cr
 \aar\quad
- \int_0^t\d s\int_0^\infty \e^{\im\tilde{M}_{s-}(a,g)-A_s(a,g)} 1_{\{H_s\le a\}} K_g(H_s,z,L_s(H_s)) \pi(\d z) \cr
 \aar\quad
+ \int_0^t\int_0^\infty \e^{\im\tilde{M}_{s-}(a,g)-A_s(a,g)} 1_{\{H_s\le a\}} K_g(H_s,z,L_s(H_s)) N(\d s,\d z) \cr
 \aar=
1 + \int_0^t\int_0^\infty \e^{\im\tilde{M}_{s-}(a,g)-A_s(a,g)} 1_{\{H_s\le a\}} \big(\e^{\im g(H_s,z,L_s(H_s))} - 1\big) \tilde{N}(\d s,\d z). \quad
 \end{align}
}The process above is a purely discontinuous square integrable $(\mcr{G}_t)$-martingale. Then
 \begin{align*}
\mbf{E}\big[\e^{\im\tilde{M}_t(a,g) - A_t(a,g)}\big]= 1.
 \end{align*}
By \eqref{intf(H_s,L_s(H_s))=} and \eqref{A_t(a,g)=def} it is easy to show that $\lim_{t\to \infty} A_t(a,g)= A(a,g)$. Then letting $t\to \infty$ in the above equality and using dominated convergence we get
 \begin{align*}
\mbf{E}\big[\e^{\im M(a,g) - A(a,g)}\big]= 1.
 \end{align*}
This gives \eqref{E[e^{iM(a,g)}]=} and completes the proof. \eproof

\bproposition\label{thW(a,f)M(a,g)ind} The random variables $W(a,f)$ and $\tilde{M}(a,g)$ given by \eqref{W(a,f)=def} and \eqref{tilM(a,g)=def} are independent and
 \begin{align}\label{Ee^{iW_t(a,f)+iM(a,g)}]}
\mbf{E}[\e^{\im W(a,f)+\im M(a,g)}]
 =
\e^{\frac{1}{2}\<W(a,f)\> + A(a,g)}, \quad a\ge 0, f\in \mcr{L}_2, g\in \mcr{L}_3(\pi).
 \end{align}
\eproposition

\bproof Clearly, the square integrable martingales defined by \eqref{e^{iW_t(a,f)+..}} and \eqref{e^{itilM_t(a,f)-..}} are orthogonal. Then we have
 \begin{align*}
\mbf{E}\big[\e^{\im W_t(a,f)+\im\tilde{M}_t(a,g)+\frac{1}{2}\<W(a,f)\>_t - A_t(a,g)}\big]= 1.
 \end{align*}
By letting $t\to \infty$ in the above equality and using dominated convergence we get
 \begin{align*}
\mbf{E}\big[\e^{\im W(a,f)+\im\tilde{M}(a,g)+\frac{1}{2}\<W(a,f)\> - A(a,g)}\big]= 1.
 \end{align*}
This proves \eqref{Ee^{iW_t(a,f)+iM(a,g)}]}, which implies the independence of $W(a,f)$ and $\tilde{M}(a,g)$ . \eproof

\medskip\noindent\textit{Proof of Proposition~\ref{th-GP-noises}.~} By Propositions~\ref{thW(a,f)=def} and~\ref{thtilM(a,g)=def}, we can define a time-space Gaussian white noise $W(\d s,\d u)$ on $(0,\infty)^2$ with intensity $\d s\d u$ and a compensated time-space Poisson random measure $\tilde{M}(\d s,\d z,\d u)$ on $(0,\infty)^3$ with intensity $\d s\pi(\d z)\d u$ such that \eqref{iintfdW=def} and \eqref{iintgdM=def} hold. It is immediate that, for $a\ge 0$, $f\in \mcr{L}_2$ and $g\in \mcr{L}_3(\pi)$,
 \begin{align*}
\int_0^a\int_0^\infty f(s,u) W(\d s,\d u)
 =
\int_0^\infty f(H_{a,s}^-,L_{a,s}^-(H_{a,s}^-))\d B_{a,s}^-
 \end{align*}
and
 \begin{align*}
\aar\int_0^a\int_0^\infty\int_0^\infty g(s,z,u) \tilde{M}(\d s,\d z,\d u) \cr
 \aar\qqquad\qqquad
= \int_0^\infty\int_0^\infty g(H_{a,s}^-,z,L_{a,s}^-(H_{a,s}^-)) \tilde{N}_a^-(\d s,\d z).
 \end{align*}
By Proposition~\ref{th(H^{a+-})-msbl}-(ii) we infer that $W(a,f)$ and $M(a,g)$ are measurable relative to $\mcr{F}_a$. Then the processes $\{W(a,f): a\ge 0\}$ and $\{M(a,g): a\ge 0\}$ are adapted to the filtration $(\mcr{F}_a)_{a\ge 0}$. For any $a,c\ge 0$, we can use \eqref{W(a,f)=def} to see that
 \begin{align}\label{iintW(a+dh,..)=}
\int_0^c\int_0^\infty f(s,u) W(a+\d s,\d u)
 \aar=
\int_0^\infty 1_{\{a<H_s\le a+c\}} f(H_s-a,L_s(H_s))\d B_s \cr
 \aar=
\int_0^\infty 1_{\{0< H_{a,s}^+-a\le c\}} f(H_{a,s}^+-a,L_{a,s}^+(H_{a,s}^+))\d B_{a,s}^+. \qquad
 \end{align}
Similarly, by \eqref{tilM(a,g)=def} we have
 \begin{align}\label{iiintM(a+dh,..)=}
\aar\int_0^c\int_0^\infty\int_0^\infty g(s,z,u) M(a+\d s,\d z,\d u) \cr
 \aar\qqquad
= \int_0^\infty\int_0^\infty 1_{\{0< H_{a,s}^+-a\le c\}} g(H_{a,s}^+-a,z,L_{a,s}^+(H_{a,s}^+)) N_a^+(\d s,\d z). \qquad
 \end{align}
The stochastic integrals in \eqref{iintW(a+dh,..)=} and \eqref{iiintM(a+dh,..)=} are $\mcr{F}^a$-measurable by Proposition~\ref{th(H^{a+-})-msbl}-(i), so they are independent of $\mcr{F}_{a+}$ by Proposition~\ref{explor-indincr}. This means that the processes $\{W(a,f): a\ge 0\}$ and $\{M(a,g): a\ge 0\}$ have independent increments relative to $(\mcr{F}_{a+})_{a\ge 0}$. Therefore $W(\d s,\d u)$ is a time-space $(\mcr{F}_{a+})$-Gaussian white noise on $(0,\infty)^2$ with intensity $\d s\d u$ and $\tilde{M}(\d s,\d z,\d u)$ is a compensated time-space $(\mcr{F}_{a+})$-Poisson random measure on $(0,\infty)^3$ with intensity $\d s\pi(\d z)\d u$. \qed

\section{The Ray--Knight theorem}

\setcounter{equation}{0}

In this section, we give the proof of the stochastic integral representation \eqref{DLEX_t(x)=x+..} for the Ray--Knight theorem. We shall identify a stochastic integral involving a $(\mcr{F}_{a+})$-adapted process with the one involving its $(\mcr{F}_{a+})$-predictable version.

%\blemma\label{th-L_{T_x}(a)adapted} \blue For any $x\ge 0$ the process $\{L_{T_x}(a): a\ge 0\}$ is $(\mcr{F}_a)$-adapted. \elemma

%\bproof \blue (A simpler proof of the result was given in \cite[Theorem~1.4.1]{DLG02}.) By Propositions~\ref{thL_{tau_a(t)}=L_{a,t}}-(ii) and~\ref{th(H^{a+-})-msbl}-(ii), the random variable $L_{a,t}^-(0)= L_{\tau_a^-(t)}(0)$ is $\mcr{F}_a$-measurable. Since $t\mapsto \tau_a^-(t)$ is strictly increasing, we see that, for any $t\ge 0$, \begin{align*}\{\eta_a^-(T_x)\ge t\}= \{T_x\ge \tau_a^-(t)\}= \{L_{\tau_a^-(t)}(0)\le x\}\in \mcr{F}_a. \end{align*} This shows that $\eta_a^-(T_x)$ is an $\mcr{F}_a$-measurable random variable. On the other hand, by Proposition~\ref{th(H^{a+-})-msbl}-(ii), the random variables $\{L_{a,t}^-(a): t\ge 0\}$ are $\mcr{F}_a$-measurable. By the continuity of $t\mapsto L_{a,t}^-(a)$ we see that $(\omega,t)\mapsto L_{a,t}^-(\omega,a)$ is a measurable map from $(\Omega\times [0,\infty),\mcr{F}_a\times \mcr{B}[0,\infty))$ to $([0,\infty), \mcr{B}[0,\infty))$. Then the composition $\omega\mapsto (\omega,\eta_a^-(\omega,T_x(\omega)))\mapsto L_{a,\eta_a^-(\omega,T_x(\omega))}^-(\omega,a)$ is $\mcr{F}_a$-measurable. By Proposition~\ref{thL_{tau_a(t)}=L_{a,t}}-(ii), the random variable $L_{T_x}(a)= L_{a,\eta_a^-(T_x)}^-(a)$ is $\mcr{F}_a$-measurable. \eproof

\blemma\label{th-iintdW=intdB} For any $x\ge 0$ and $a\ge 0$ we have
 \begin{align}\label{iintdW=intdB}
\int_0^a\int_0^{L_{T_x}(v)} W(\d v,\d u)
 =
\int_0^\infty 1_{\{H_s\le a\}}1_{\{L_s(H_s)\le L_{T_x}(H_s)\}} \d B_s.
 \end{align}
\elemma

\bproof Fix $a> 0$ and write $a_k= a/2^k$ for $k\ge 1$. For $v\ge 0$ set $q_k(v)= \max([0,v)\cap \{ia_k: i=0,1,2,\cdots\})$ with $\max(\emptyset)= 0$ by convention. Then $q_k(v)\to v$ increasingly as $k\to \infty$. Since $L_{T_x}(ia_k)$ is independence of $W(ia_k+\d v,\d u)$, we can apply \eqref{iintW(a+dh,..)=} to the function $f(s,u)= 1_{\{u\le L_{T_x}(ia_k)\}}$ to see that
 \begin{align}\label{iintl_kdW=intldB}
\int_0^a\int_0^{L_{T_x}(q_k(v))} W(\d v,\d u)
 \aar=
\sum_{i=0}^{2^k-1}\int_0^{a_k}\int_0^{L_{T_x}(ia_k)} W(ia_k+\d v,\d u) \cr
 \aar=
\sum_{i=0}^{2^k-1}\int_0^\infty 1_{\{H_{ia_k,s}^+-ia_k\le a_k\}} 1_{\{L_{ia_k,s}^+(H_{ia_k,s}^+)\le L_{T_x}(ia_k)\}} \d B_{ia_k,s}^+ \cr
 \aar=
\sum_{i=0}^{2^k-1}\int_0^\infty 1_{\{ia_k< H_s\le (i+1)a_k\}} 1_{\{L_s(H_s)\le L_{T_x}(ia_k)\}} \d B_s \cr
 \aar=
\int_0^\infty 1_{\{H_s\le a\}}1_{\{L_s(H_s)\le L_{T_x}(q_k(H_s))\}} \d B_s.
 \end{align}
Moreover, we have
 \begin{align*}
\aar\mbf{E}\Big[\Big|\int_0^a\int_0^{L_{T_x}(v)} W(\d v,\d u) - \int_0^t\int_0^{L_{T_x}(q_k(v))} W(\d v,\d u)\Big|^2\Big] \cr
 \aar\qquad
= \mbf{E}\Big[\int_0^a\d v\int_0^\infty\big|1_{(0,L_{T_x}(v)]}(u) - 1_{(0,L_{T_x}(q_k(v))]}(u)\big| \d u\Big] \cr
 \aar\qquad
= \mbf{E}\Big[\int_0^a|L_{T_x}(v)-L_{T_x}(q_k(v))|\d v\Big]
 \end{align*}
and
 \begin{align*}
\aar\mbf{E}\Big[\Big|\int_0^\infty 1_{\{H_s\le a\}}\big(|1_{\{L_s(H_s)\le L_{T_x}(H_s)\}} - 1_{\{L_s(H_s)\le L_{T_x}(q_k(H_s))\}}\big) \d B_s\Big|^2\Big] \cr
 \aar\qquad
= \mbf{E}\Big[\int_0^\infty 1_{\{H_s\le a\}}\big|1_{\{L_s(H_s)\le L_{T_x}(H_s)\}} - 1_{\{L_s(H_s)\le L_{T_x}(q_k(H_s))\}}\big|\d s\Big] \cr
 \aar\qquad
= \mbf{E}\Big[\int_0^a \d v\int_0^\infty \big|1_{\{L_s(v)\le L_{T_x}(v)\}} - 1_{\{L_s(v)\le L_{T_x}(q_k(v))\}}\big|\d L_s(v)\Big] \cr
 \aar\qquad
= \mbf{E}\Big[\int_0^a \d v\int_0^\infty \big|1_{\{u\le L_{T_x}(v)\}} - 1_{\{u\le L_{T_x}(q_k(v))\}}\big|\d u\Big] \cr
 \aar\qquad
= \mbf{E}\Big[\int_0^a |L_{T_x}(v)-L_{T_x}(q_k(v))|\d v\Big].
 \end{align*}
Since the CB-process $\{L_{T_x}(a): a\ge 0\}$ has at most countably many jumps, by \eqref{E[L_{T_x}(a)]=xe^{}} and dominated convergence it is easy to see that
 \begin{align*}
\mbf{E}\Big[\int_0^a |L_{T_x}(v)-L_{T_x}(q_k(v))|\d v\Big]\to 0, \quad k\to \infty.
 \end{align*}
Then taking the $L^2(\mbf{P})$-limits as $k\to \infty$ in \eqref{iintl_kdW=intldB} we obtain \eqref{iintdW=intdB}. \eproof

\blemma\label{th-iiintzdM=iintzdN} For any $x\ge 0$ and $a\ge 0$ we have
 \begin{align}\label{iiintzdM=iintzdN}
\int_0^a\int_0^\infty\int_0^{L_{T_x}(v)} z \tilde{M}(\d v,\d z,\d u)
 =
\int_0^\infty\int_0^\infty 1_{\{H_s\le a\}}1_{\{L_s(H_s)\le L_{T_x}(H_s)\}}z \tilde{N}(\d s,\d z).
 \end{align}
\elemma

\bproof Suppose that the Poisson random measure $N(\d s,\d z)$ is given by \eqref{N=sum_I}. In view of \eqref{iintgdM=def}, we have
 \begin{align}\label{M(ds,dz,du)=rep}
M(\d s,\d z,\d u)= \sum_{i\in I} \delta_{(H_{t_i},z_i,L_{t_i}(H_{t_i}))}(\d s,\d z,\d u).
 \end{align}
Then, for any $A\ge \delta> 0$,
 \begin{align*}
\int_0^a\int_\delta^A\int_0^{L_{T_x}(v)} z M(\d v,\d z,\d u)
 \aar=
\int_0^\infty\int_\delta^A\int_0^\infty 1_{\{v\le a\}} z 1_{\{u\le L_{T_x}(v)\}} M(\d v,\d z,\d u) \cr
 \aar=
\int_0^\infty\int_\delta^A 1_{\{H_s\le a\}} 1_{\{L_s(H_s)\le L_{T_x}(H_s)\}}z N(\d s,\d z).
 \end{align*}
It follows that
 \begin{align}\label{int_del^AzdtM=..}
\int_0^a\int_\delta^A\int_0^{L_{T_x}(v)} z \tilde{M}(\d v,\d z,\d u)
 =
\int_0^\infty\int_\delta^A 1_{\{H_s\le a\}}1_{\{L_s(H_s)\le L_{T_x}(H_s)\}}z \tilde{N}(\d s,\d z).
 \end{align}
Now, letting $\delta= 1$ and taking the $L^1(\mbf{P})$-limits as $A\to \infty$ in \eqref{int_del^AzdtM=..}, we have
 \begin{align*}
\int_0^a\int_1^\infty\int_0^{L_{T_x}(v)} z \tilde{M}(\d v,\d z,\d u)
 =
\int_0^\infty\int_1^\infty 1_{\{H_s\le a\}}1_{\{L_s(H_s)\le L_{T_x}(H_s)\}}z \tilde{N}(\d s,\d z).
 \end{align*}
Similarly, letting $A= 1$ and taking the $L^2(\mbf{P})$-limits as $\delta\to 0$ in \eqref{int_del^AzdtM=..}, we get
 \begin{align*}
\int_0^a\int_0^1\int_0^{L_{T_x}(v)} z \tilde{M}(\d v,\d z,\d u)
 =
\int_0^\infty\int_0^1 1_{\{H_s\le a\}}1_{\{L_s(H_s)\le L_{T_x}(H_s)\}}z \tilde{N}(\d s,\d z).
 \end{align*}
Then the equality \eqref{iiintzdM=iintzdN} holds. \eproof

\medskip\noindent\textit{Proof of Theorem~\ref{th-DL-inteq}.~} By Lemma~\ref{th1_{sleT_x}=a.s.}, the left-continuous process $s\mapsto 1_{\{s\le T_x\}}$ is a $(\mcr{G}_s)$-predictable version of $s\mapsto 1_{\{L_s(H_s)\le L_{T_x}(H_s)\}}$. It follows that
 \begin{align*}%\label{L_{T_x}(a)=x+..}
\int_0^{T_x} 1_{\{H_s\le a\}}\d \xi_s\aar= x - \alpha\int_0^\infty 1_{\{H_s\le a, L_s(H_s)\le L_{T_x}(H_s)\}} \d s \cr
 \aar\quad
+ \sqrt{2\beta}\int_0^\infty 1_{\{H_s\le a, L_s(H_s)\le L_{T_x}(H_s)\}}\d B_s \cr
 \aar\quad
+ \int_0^\infty\int_0^\infty 1_{\{H_s\le a, L_s(H_s)\le L_{T_x}(H_s)\}} z \tilde{N}(\d s,\d z).
 \end{align*}
By Lemmas~\ref{th-iintdW=intdB} and~\ref{th-iiintzdM=iintzdN}, we can rewrite the above equation as
 \begin{align*}
L_{T_x}(a)\aar= x - \alpha\int_0^a L_{T_x}(v) \d v + \sqrt{2\beta}\int_0^a\int_0^{L_{T_x}(v)} W(\d v,\d u) \cr
 \aar\quad
+ \int_0^a\int_0^\infty\int_0^{L_{T_x}(v)} z \tilde{M}(\d v,\d z,\d u).
 \end{align*}
This $\{L_{T_x}(a): a\ge 0\}$ solves the stochastic equation \eqref{DLEX_t(x)=x+..}. The pathwise uniqueness of the solution holds by \cite[Theorem~3.1]{DaL12}. \qed

\bigskip

\textbf{Acknowledgements.~} We are grateful to Professor Elie A\"id\'ekon, who pointed out an error in an earlier form of Lemma~\ref{th-iintdW=intdB}. We thank Professors Jean Bertoin and Thomas Duquesne for their helpful comments on Proposition~\ref{th-betL_t(0,R)=..}. This research is supported by the National Key R{\&}D Program of China (No.~2020YFA0712901) and the Natural Science Foundation of China (No.~12271029 and No.~12571153).

\end{document}